\documentclass{ifacconf}

\usepackage{graphicx}      % include this line if your document contains figures
\usepackage{natbib}        % required for bibliography
\usepackage{algorithm}  %algorithm
\usepackage{mathtools,color}	%coloneqq
\usepackage{amsmath,amssymb}

\usepackage{siunitx}  
\usepackage{pgfplots} %matlab2tikz

\pgfplotsset{compat=newest}
\usetikzlibrary{plotmarks}
\usetikzlibrary{arrows.meta}
\usepgfplotslibrary{patchplots}
\usepackage{grffile}

\DeclareMathOperator*{\argmin}{argmin} % no space, limits underneath in displays

\newtheorem{assu}{Assumption}
\newtheorem{rema}{Remark}
\begin{document}
\begin{frontmatter}

\title{A Distributed Active Set Method for Model Predictive Control} 

\author[First]{Gösta Stomberg} 
\author[First]{Alexander Engelmann} 
\author[First]{Timm Faulwasser}

\address[First]{Institute for Energy Systems, Energy Efficiency and Energy Economics, Faculty of Electrical Engineering and Information Technology, TU Dortmund University, 44227 Dortmund, Germany (e-mail: \{goesta.stomberg,alexander.engelmann\}@tu-dortmund.de, timm.faulwasser@ieee.org)}

\begin{abstract}                % Abstract of not more than 250 words.
This paper presents a novel distributed active set method for  model predictive control of linear systems. The method combines a primal active set strategy with a decentralized conjugate gradient method to solve convex quadratic programs. An advantage of the proposed method compared to existing distributed model predictive algorithms is the primal feasibility of the iterates. Numerical results show that the proposed method can compete with the alternating direction method of multipliers in terms of communication requirements for a chain of masses example.
\end{abstract}

\begin{keyword}
distributed control, model predictive control, active set methods, conjugate gradient methods, ADMM
\end{keyword}

\end{frontmatter}
%===============================================================================

\section{Introduction}
In many control applications such as process control, distributed model predictive control (DMPC) is promising due to the absence of central coordination and a limited information exchange.
DMPC schemes are often realized by solving  optimal control problems (OCPs)---respectively, appropriate discretizations thereof---by means of distributed optimization algorithms.
A challenge is that many distributed optimization algorithms guarantee consensus constraint satisfaction  only asymptotically. 
This implies that per sampling instant usually many distributed optimization iterations have to be executed in order to ensure closed-loop properties such as stability and recursive feasibility \citep{Kogel2012,Mota2012,Pu2014,Rostami2017}. 
To overcome this limitation and to enable early termination of the optimization, constraint-tightening approaches have been proposed in \cite{Doan2011,Giselsson2014}.
\cite{Wang2017} extend these approaches to state constraints. However, constraint tightening potentially leads to a loss in control performance, or to a reduced domain of attraction resulting from the tightened constraint set.

A second line of research enables decomposition via coordinate descent methods, where each subsystem solves its own OCP while keeping the influence of the neighboring systems fixed. 
Such methods guarantee feasible iterates \citep{Necoara2012}.
A  Jacobi-type method has been proposed by \cite{Stewart2010}, where coupling variables are averaged between subsystems and
a linear rate of convergence has been shown by \cite{Gross2013}.
This approach was extended to event-based communication in \cite{Gross2016}.
A condition for verifying convergence a posteriori is presented in 
\citep{DangDoan2017}, where the authors also provide convergence guarantees for the special case of chain-linked systems.
However, a general shortcoming of Jacobi-type iterations is that---although they decrease the objective in every iteration---convergence to a minimizer is often not guaranteed.
Moreover, for fast convergence they require carefully chosen weights. 
For a  comprehensive overview  on  DMPC approaches  we refer to \cite{Muller2017}.

This paper aims to overcome the conservatism introduced by constraint tightening.
We introduce a novel distributed active set method (ASM) that builds upon a decentralized conjugate gradient~(DCG) method. The method converges to the \emph{exact solution} of an OCP in a finite number of iterations and thus avoids constraint tightening. Moreover, the proposed method has the advantage that the iterates remain feasible after a feasible initial point has been computed. This makes early termination possible without losing recursive feasibility guarantees. Furthermore,  our method is tuning-free, has a practically super-linear convergence rate and numerical examples indicate a low communication footprint.

The paper is structured as follows. Section \ref{sec:ProblemStatement} formulates the DMPC scheme for the control of a linear system network. Section \ref{sec:ProposedMethod} presents the proposed distributed ASM for solving the arising QPs in the DMPC scheme. Section \ref{sec:Results} presents numerical results on an example chain of masses system where our method is compared to ADMM. Appendix \ref{app:ADMM} briefly discusses the ADMM variant we use as a benchmark.

\section{Problem Statement} \label{sec:ProblemStatement}

We consider a network of linear time-invariant systems $i \in \mathcal{M}=\{1,\dots,M\}$ with dynamics
\begin{equation}
x_i^{+}=A_{ii}x_i+B_iu_i + \sum_{j\in \mathcal{M}_i^{\text{in}}}A_{ij}x_j, \quad x_i(0) = x_{i,0}.\label{eq:dyn}\\ 
\end{equation} where $x_i\in \mathbb{R}^{n_i}$ and $u_i \in \mathcal{U}_i \subset \mathbb{R}^{m_i}$ denote the state and input of agent $i$. We use a subscript $(\cdot)_i$ throughout the paper to mark a variable's association to agent $i$. The matrices $A_{ij}\in \mathbb{R}^{n_i\times n_j}$ create state couplings between neighboring agents. We distinguish between in-neighbors $\mathcal{M}_i^{\text{in}} = \{j \in \mathcal{M}\,|\, A_{ij} \neq 0\}$ and out-neighbors $\mathcal{M}_i^{\text{out}} = \{j\in \mathcal{M}\,|\,A_{ji} \neq 0\}$ and we denote their union by $\mathcal{M}_i \doteq \mathcal{M}_i^{\text{in}} \cup \mathcal{M}_i^{\text{out}}$. The dynamics \eqref{eq:dyn} can be rewritten as 
\begin{equation} \label{eq:distDyn}
x_i^{+}=A_{ii}x_i+B_iu_i + \sum_{j\in \mathcal{M}_i^{\text{in}}}A_{ij}v_{ji}, \quad x_i(0) = x_{i,0}
\end{equation} where $v_{ji} \doteq x_j$ are state copies. We gather the state copies of agent $i$'s in-neighbors in
\begin{equation*}
v_i \doteq \text{col}_{j \in \mathcal{M}_i^{\text{in}}} (v_{ji}) \in \mathbb{R}^{v_i}.
\end{equation*}

\begin{assu}[Input Constraints]
	\label{assu:constraintSet}
	The sets $\mathcal{U}_i$ are compact convex polytopes and contain the origin in their interior. \hfill $\square$
\end{assu}

\begin{rema}[State Constraints]
	For the sake of simplicity, we consider only input constraints in the present paper. Our approach can be extended to include state constraints as well as terminal constraints. \hfill $\square$
\end{rema}

Let $x\doteq \text{col}_{i\in \mathcal{M}}(x_i) \in \mathbb{R}^n$ denote the overall system state and $u\doteq \text{col}_{i\in \mathcal{M}}(u_i) \in \mathbb{R}^m$ the overall input. The overall system dynamics are $x^+ = A x + B u$, where $A\in \mathbb{R}^{n\times n}$ is a block matrix with entries $A_{ij}$ and $B = \text{diag}_{i\in \mathcal{M}}(B_i) \in \mathbb{R}^{n \times m}$. 

\begin{assu}[Controllability]
	\label{assu:ctrb}
	The overall system $(A,B)$ is controllable. \hfill $\square$
\end{assu}

\subsection{MPC Formulation}
Let $\boldsymbol{x}_i$, $\boldsymbol{u}_i$ and $\boldsymbol{v}_i$ denote the state, input and copy sequences $\boldsymbol{x}_i^\top \doteq  [x_i^{0\top},\dots, x_i^{N-1\top}] $, $\boldsymbol{u}_i^\top \doteq  [u_i^{0\top},\dots, u_i^{N-1\top}] $ and $\boldsymbol{v}_i^\top \doteq  [v_i^{0\top},\dots, v_i^{N-1\top}] $ for agent $i$ with horizon $N$. The OCP reads
\begin{subequations}
	\label{eq:centralOCP}
\begin{align}
&\min_{\boldsymbol{x}_i, x_i^N, \boldsymbol{u}_i, \boldsymbol{v}_i} \quad \sum_{i \in \mathcal{M}} J_i(\boldsymbol{x}_i,x_i^N,\boldsymbol{u}_i)\\
\notag \textrm{su}& \textrm{bject} \textrm{ to}  \textrm{ for all } i \in \mathcal{M}:\\
 \quad x_i^{+}&=A_{ii}x_i+B_iu_i + \sum_{j\in \mathcal{M}_i^{\text{in}}}A_{ij}v_{ji}, \; x_i^0 = x_{i,0},\\
u_i &\in \mathcal{U}_i, \label{eq:inCnstr}\\
v_{ji} &= x_j \quad \forall j \in \mathcal{M}_i^{\text{in}},\label{eq:couplingOCP}
\end{align}
\end{subequations} where the objective for agent $i \in \mathcal{M}$ is given by
\begin{equation*}
J_i(\boldsymbol{x}_i, x_i^N, \boldsymbol{u}_i) \doteq \frac{1}{2}x_i^{N\top}P_i x_i^N+\frac{1}{2} \sum_{k=1}^{N-1}(x_i^{k\top}Q_ix_i^k + u_i^{k\top}R_iu_i^k).
\end{equation*} 
\begin{assu}
	\label{assu:weights}
	The weighting matrices $Q_i$ and $R_i$ are positive definite. The matrices $P_i$ are positive semi-definite. \hfill $\square$
\end{assu} 

Given measured or observed states $x_{i,0}$, the MPC scheme solves this OCP for the optimal input sequences $\boldsymbol{u}_i^\ast$ and then applies the first part of the input sequence $u_i^{0,\ast}$ as a control input to each agent. This MPC scheme asymptotically stabilizes the origin given Assumptions \ref{assu:constraintSet}--\ref{assu:weights} for sufficiently long horizons $N$ on a set of initial conditions $\mathbb{X}_0$. The set $\mathbb{X}_0$ on which closed-loop stability can be guaranteed depends on the interplay of the input constraints $\mathcal{U}$ and the stability of $A$, cf. \citep{Boccia_2014}.

\subsection{Distributed OCP Formulation}
We next aggregate the state, input and copy trajectories into decision variables for each agent into $z_i^\top \doteq [\boldsymbol{x}_i^\top,x_i^{N\top},\boldsymbol{u}_i^\top,\boldsymbol{v}_i^\top] \in \mathbb{R}^{N_{zi}}$ where $\boldsymbol{v}_i^\top \doteq [v_i^{0\top},\dots,v_i^{N-1\top}] \in \mathbb{R}^{Nv_i}$ and $n_{zi} \doteq (N+1)n_i + N (m_i + v_i)$. The decision variables $z_i$ are not typeset bold to simplify the notation even though they contain trajectories. With this we rewrite OCP \eqref{eq:centralOCP} as the partially separable quadratic program~(QP)
\begin{subequations}
	\label{ocp:distributed}
\begin{align}
\label{eq:OCP}
\min_{z_i} \quad  \sum_{i\in \mathcal{M}} &\frac{1}{2}z_i^{\top}H_i z_i\\
\textrm{s.t.} \quad C_i^{\mathcal{E}} z_i &= b_i^{\mathcal{E}} \quad \forall i \in \mathcal{M}, \label{eq:eqCon}\\
C_i^{\mathcal{I}} z_i &\leq b_i^{\mathcal{I}}\quad \forall i \in \mathcal{M},\label{eq:IneqCon}\\
 \sum_{i\in \mathcal{M}} C_i^c z_i &= 0. \label{eq:cCon}
\end{align}
\end{subequations} The Hessian  matrix is given by 
\begin{equation}
H_i\doteq\text{diag}(\boldsymbol{\tilde{Q}}_i,P_i,\boldsymbol{R}_i,[\boldsymbol{\tilde{Q}}_j]_{j\in \mathcal{M}_i^{\text{in}}}) \in \mathbb{R}^{n_{zi} \times n_{zi}},
\end{equation} where $\boldsymbol{\tilde{Q}}_i \doteq 1/(\lvert \mathcal{M}_i^{\text{out}}\rvert + 1) \cdot \text{diag}(Q_i,\dots,Q_i) \in \mathbb{R}^{Nn_i \times Nn_i}$ and $\boldsymbol{R}_i \doteq \text{diag}(R_i,\dots,R_i) \in \mathbb{R}^{Nm_i \times Nm_i}$. 
This choice of $\boldsymbol{\tilde{Q}}_i$ spreads the cost associated with state $x_i$ evenly among agent $i$ and its outgoing neighbors. The equality constraints~\eqref{eq:eqCon} with $C_i^{\mathcal{E}} \in \mathbb{R}^{n_{hi}\times n_{zi}}$ include the initial condition  and dynamics \eqref{eq:distDyn}; the inequality constraints \eqref{eq:IneqCon} include the input constraints \eqref{eq:inCnstr}; and the coupling constraints (\ref{eq:cCon}) with $C_i^c \in \mathbb{R}^{n_c \times n_{zi}}$ include the coupling between states and copies \eqref{eq:couplingOCP}. We propose to solve OCP (\ref{ocp:distributed}) with the following distributed ASM.

\section{distributed active set method}\label{sec:ProposedMethod}
Active set methods (ASMs) are well-known for solving inequality-constrained convex QPs. 
In particular, primal ASMs produce primal feasible iterates \citep{Nocedal2006}, which is advantageous for MPC. We next present a distributed primal ASM for solving OCP (\ref{ocp:distributed}).

\subsection{Distributed Active Set Method}
Let $c_i^{j\mathcal{I}}$ denote the $j^\text{th}$ row of $C_i^{\mathcal{I}}$ and $b_i^{j\mathcal{I}}$ the $j^\text{th}$ element of $b_i^{\mathcal{I}}$.  The active set of agent $i$ is given by
\begin{equation*} 
	\mathcal{A}(z_i) \doteq \{j \,|\, c_i^{j\mathcal{I}} z_i = b_i^{j\mathcal{I}}\}.
\end{equation*}

The proposed method first chooses an initial active set $\mathcal{A}(z_i^0)$ and a corresponding iterate $z_i^0$ for each agent. 
Afterwards, the method takes steps $z_i^{n+1}\doteq z_i^n + \alpha^n \Delta z_i^n$ until it has converged. 
The step directions $\Delta z_i^n$ are obtained by solving the equality constrained QP
\begin{subequations} \label{eq:EQP}
\begin{align}
\min_{\Delta z_i} \quad  \sum_{i\in \mathcal{M}} \frac{1}{2}(\Delta z_i^{\top}&H_i \Delta z_i + g_i^{n\top} \Delta z_i)\\
\textrm{s.t.} \quad C_i^n \Delta z_i &= d_i && \hspace{-2.5cm}| \;\gamma_i \quad \forall i \in \mathcal{M},\label{eq:EQPec}\\
\sum_{i\in \mathcal{M}} C_i^c \Delta z_i &= 0 && \hspace{-2.5cm}| \; \lambda_C \label{eq:EQPcc}
\end{align}
\end{subequations} where $C_i^n \doteq \begin{bmatrix}
C_i^{\mathcal{E}}\\
\text{col}_{j\in \mathcal{A}(z_i)}[c_i^{j\mathcal{I}}]\end{bmatrix}$, $g_i^n \doteq H_i z_i^n$ and $d_i=0$. 
Here, $\gamma_i$ and $\lambda_C$ are the Lagrange multipliers associated with the constraints (\ref{eq:EQPec}) and (\ref{eq:EQPcc}). Then, each agent computes the largest step length $\alpha_i^n \in (0,1]$ such that primal feasibility is maintained \citep{Nocedal2006}:
\begin{equation} \label{eq:alpha}
\alpha_i^n \doteq \min\,\left\{1, \min_{j \notin \mathcal{A}(z_i^n), c_i^{j\mathcal{I}}p_i^n<0} \frac{b_i^{j\mathcal{I}} - c_i^{j\mathcal{I}} z_i^n}{c_i^{j\mathcal{I}}p_i^n}\right\}.
\end{equation} The step length $\alpha^n = \min \{\alpha_1^n,\dots,\alpha_M^n\}$ is then chosen for all agents to obtain $z_i^{n+1}$. This choice ensures that $z_i^{n+1}$ remains feasible for all agents. If $\alpha<1$, then an inactive constraint of one agent is blocking and the constraint is added to the respective agent's active set. This process is repeated until $\lVert \Delta z_i^n \rVert < \varepsilon \quad \forall i \in \mathcal{M}$ where $\varepsilon$ is small. 
Then, each agent checks for dual feasibility
\begin{equation*}
\gamma_i^{j\mathcal{I}} \geq 0, \quad \forall j \in \mathcal{A}(z_i^n).
\end{equation*} If dual feasibility is attained for all agents, then $z^n$ is returned as solution to problem (\ref{eq:EQP}). Else the constraint corresponding to the smallest Lagrange multiplier among all agents is removed from the respective agent's active set. Algorithm 1 summarizes the ASM and it further contains the condensing and DCG steps that we subsequently describe in Subsections \ref{sec:condensing} and \ref{sec:DCG}. 

\begin{rema}[ASM initialization]\label{rema:initASM}
Different ways of initializing ASMs with $\mathcal{A}(z_i^0)$ and $z_i^0$ have been reported in the literature. \cite{Ferreau2008} present an ASM that is initialized on a homotopy path between the current and the previous sample. \cite{Klauco2019} propose to warm-start ASMs with classification methods from supervised learning. Here, we initialize with the optimal active set  from the previous MPC iteration as a warm-start or with $\mathcal{A}(z_i^0)=\emptyset$ if no previous MPC iteration is available. We then find $z_i^0$ by solving a QP similar to \eqref{ocp:distributed} but with \eqref{eq:IneqCon} replaced by the active set. If the obtained $z_i^0$ violates an inequality constraint, then this constraint is added to the active set. This is repeated until a feasible initialization is obtained. \hfill $\square$
\end{rema} 

\begin{rema}[Early termination]\label{rema:earlyTermination}
	The active set method produces primal feasible iterates that satisfy \eqref{eq:eqCon}--\eqref{eq:cCon}. It can therefore be terminated early and still guarantee stability for MPC schemes with terminal constraints where feasibility implies stability \citep{Scokaert1999}. \hfill $\square$
\end{rema}

\begin{algorithm}[h]
	%\small 
	\caption[fontsize=Large]{Distributed ASM}
	\textbf{Initialization: $ \mathcal{A}(z_i^0)$ and $z_i^0$.}\\
	\textbf{Repeat until convergence:}\\	
	\hspace*{1em} Condense QP (\ref{eq:EQP}) into \eqref{eq:Schur}.\\
	\hspace*{1em} Solve \eqref{eq:Schur} with DCG for $\lambda_C^n$.\\
	\hspace*{1em} Obtain $\Delta z_i$ via backsubstitution \eqref{eqn:backsubstitution} for all $i\in \mathcal{M}$.\\
	\hspace*{1em} \textbf{If} $\lVert \Delta z_i \rVert < \varepsilon \quad \forall i\in \mathcal{M}$:\\
	\hspace*{2em} $\gamma_i^n = (C_i^n C_i^{n\top})^{-1}C_i^n(-g_i-C_i^{c\top} \lambda_C^n) \,\, \forall i\in \mathcal{M}$.\\
	\hspace*{2em} \textbf{If} $\gamma_i^{j\mathcal{I}} \geq 0 \quad \forall j \in \mathcal{A}(z_i^n), \quad \forall i\in \mathcal{M}$:\\
	\hspace*{3em} \textbf{Return $z_i^n$ and terminate}.\\
	\hspace*{2em} \textbf{Else}:\\
	\hspace*{3em} Find the smallest $\gamma_i^{j\mathcal{I}}$ among all $i\in \mathcal{M}$.\\
	\hspace*{3em} Remove the respective constraint $j$ from $\mathcal{A}(z_i^n).$\\
	\hspace*{1em} \textbf{Else}:\\
	\hspace*{2em} Compute $\alpha_i^n$ locally according to (\ref{eq:alpha}) for all $ i\in \mathcal{M}$.\\
	\hspace*{2em} Determine agent with the smallest $\alpha_i$ and set $\alpha=\alpha_i$.\\
	\hspace*{2em} $\mathcal{A}(z_i^n) \leftarrow \mathcal{A}(z_i^n) \cup \{j\}$ where $j$ is the blocking\\ \hspace*{2em} constraint.\\
	\hspace*{2em} Local step: $z_i^{n+1} \leftarrow z_i^n + \alpha^n \Delta z_i^n$ for all $i\in \mathcal{M}$.\\
	\hspace*{1em} \textbf{End If} 
\end{algorithm}

\subsection{Condensing the QP} \label{sec:condensing}
The core idea of the distributed ASM is to solve QP (\ref{eq:EQP}) with a decentralized conjugate gradient method (DCG) from \citep{Engelmann2021}. To apply DCG, we condense and rewrite QP~(\ref{eq:EQP}) as
\begin{equation}
\label{eq:Schur}
\left(\sum_{i\in \mathcal{M}} S_i \right) \lambda_C = \sum_{i\in \mathcal{M}} s_i,
\end{equation} where $S_i$ and $s_i$ are yet to be defined. In particular, $\sum_{i\in \mathcal{M}} S_i$ is positive definite and the matrices $S_i$ have favorable sparsity properties that motivate DCG.

To arrive at \eqref{eq:Schur}, we first apply the nullpace method \citep{Nocedal2006} to eliminate the equality constraints $C_i \in \mathbb{R}^{n_{hi}\times n_{zi}}$. Let 
\begin{equation}
\label{eq:null}
\Delta z_i = Z_i v_i + Y_i w_i,
\end{equation} where the columns of $Z_i \in \mathbb{R}^{n_{zi}\times (n_{zi}-n_{hi})}$ form a nullspace of $C_i$ and $Y_i \in \mathbb{R}^{n_{zi}\times n_{hi}}$ is chosen such that $\begin{bmatrix}Z_i&Y_i\end{bmatrix}$ is invertible. This choice of $Y_i$ together with $C_iZ_i=0$ result in $C_i Y_i$ being nonsingular and inserting into the equality constraints yields $w_i = (C_i Y_i)^{-1}d_i$. We insert (\ref{eq:null}) into (\ref{eq:EQP}) and obtain
\begin{subequations}\label{eq:EQP3}
\begin{align}
\min_{v_i} \quad  & \sum_{i\in \mathcal{M}} (\frac{1}{2}v_i^{\top}\bar{H}_i v_i + \bar{g}_i^{\top} v_i)\\
\textrm{s.t.} \;\;& \sum_{i\in \mathcal{M}} (\bar{C}_i^c v_i + b_i) = 0 \;\;\mid \lambda_C
\end{align}
\end{subequations} where $\bar{H}_i \doteq Z_i^{\top}H_i Z_i$, $\bar{g}_i \doteq Z_i^{\top}g_i + Z_i^\top H_i Y_i w_i$, $\bar{C}_i^c \doteq C_i^c Z_i$ and $b_i\doteq C_i^c Y_iw_i$. The KKT conditions of (\ref{eq:EQP3}) read
\begin{equation}
\label{eq:KKT}
\begin{bmatrix}
	\bar{H}_1 & & & \bar{C}_1^{c\top}\\
			  & \ddots & & \vdots\\
			  & & \bar{H}_M & \bar{C}_M^{c\top}\\
	\bar{C}_1^c & \dots & \bar{C}_M^c & 0		  
\end{bmatrix} \begin{bmatrix}
v_1\\
\vdots\\
v_M\\
\lambda_C
\end{bmatrix} = \begin{bmatrix}
-\bar{g}_1\\
\vdots\\
-\bar{g}_M\\
-\textstyle \sum_i b_i
\end{bmatrix}.
\end{equation} 

\begin{assu}
	The matrices $\bar{H}_i$ are positive definite. \hfill $\square$
\end{assu}

Utilizing that $\bar{H}_i$ is positive definite 
for all $i \in \mathcal{M}$, we solve~(\ref{eq:KKT}) for $v_i = \bar{H}_i^{-1}(-\bar{g}_i-\bar{C}_i^{c\top} \lambda_C)$. Inserting back into the KKT conditions yields (\ref{eq:Schur}) where $S_i \doteq \bar{C}_i^c \bar{H}_i^{-1} \bar{C}_i^{c\top}$ and $s_i \doteq b_i - \bar{C}_i^c \bar{H}_i^{-1} \bar{g}_i.$ The QP (\ref{eq:EQP}) has hence been reduced to a set of $n_c$ linear equations. The solution to QP~\eqref{eq:EQP} can be obtained via backsubstitution: 
\begin{equation}\label{eqn:backsubstitution}
\Delta z_i = Z_i\bar{H}_i^{-1}(-\bar{g}_i-\bar{C}_i^{c\top} \lambda_C) + Y_i(C_i Y_i)^{-1}d_i. 
\end{equation}

\begin{rema}\label{rema:d}
	The above derivation allows for $d_i \neq 0$ in \eqref{eq:EQPec}. This is needed to perform the ASM initialization with DCG as described in Remark \ref{rema:initASM}. Once the ASM initialization is complete, we use DCG to find $\Delta z_i$ and in this case $d_i = 0$. \hfill $\square$
\end{rema}

\subsection{The Decentralized Conjugate Gradient Method}\label{sec:DCG}
Next, we recall the decentralized conjugate gradient (DCG) method from \citep{Engelmann2021} to solve the positive definite linear system of equations (\ref{eq:Schur}) in decentralized fashion. Consider the system of equations 
\begin{equation*}
S \lambda_C = s,
\end{equation*} where $S \doteq \sum_i S_i$ is positive definite and $s \doteq \sum_i s_i$. The centralized conjugate gradient method \citep{Nocedal2006} solves this system iteratively for $\lambda_C$ until the residual at iteration $n$ defined as $r^n \doteq s-S\lambda_C^n$ vanishes. The iterations read
\begin{subequations}
\begin{align}
\alpha^n &= \frac{r^{n\top}r^n}{p^{n\top}Sp^n},\label{eq:CGalpha}\\
\lambda_C^{n+1} &= \lambda_C^n + \alpha^n p^n,\label{eq:CGlambda}\\
r^{n+1} &= r^n - \alpha^n S p^n,\label{eq:CGr}\\
\beta^n &= \frac{r^{n+1\top}r^{n+1}}{r^{n\top}r^n},\\
p^{n+1} &= r^{n+1} + \beta^n p^n\label{eq:CGp},
\end{align}
\end{subequations} with the initialization $r^0 = p^0 = s-S \lambda_C^0$. The idea of DCG is to introduce local versions of the CG variables $\lambda_C, r$ and $p$ and to decompose the CG updates into local updates. Since DCG is a reformulation of the centralized CG method, DCG is guaranteed to solve the QP \eqref{eq:EQP} in a finite number of iterations and exhibits a practically superlinear convergence rate.

We exploit sparsity properties of the $S_i$ matrices that result from zero rows in $C_i^c$. The rows of $C_i^c$ couple variables between two agents: the states of one agent and copies of those states in a neighboring agent. Each row in $C_i^c$ represents a coupling constraint and the row is non-zero if and only if agent $i$ is coupled via the associated constraint. Observe that zero rows of $C_i^c$ lead to zero rows and columns in $S_i$, because $S_i = C_i^c Z_i \bar{H}_i Z_i^{\top}C_i^{c\top}$. Hence, only those rows and columns of $S_i$ are non-zero that correspond to the coupling constraints that couple agent $i$ to its neighbors $\mathcal{M}_i.$ Let $c_i^{j,c}$ denote the $j^\text{th}$ row of $C_i^c$. We let ${\mathcal{C}(i)\doteq \{j\,|\,c_i^{j,c}\neq 0\}}$ denote the consensus constraints that couple agent $i$ to its neighbors. We next introduce matrices $I_{\mathcal{C}(i)} \in \mathbb{R}^{\lvert \mathcal{C}(i) \rvert \times n_c}$ that map from global CG variables to local variables. The matrices $I_{\mathcal{C}(i)}$ are obtained by taking the identity matrix $I \in \mathbb{R}^{n_c \times n_c}$ and by subsequently eliminating all rows that do not belong to any constraint in $\mathcal{C}(i)$. With this we introduce the local variables
\begin{align*}
	\lambda_{C,i} &\doteq I_{\mathcal{C}(i)} \lambda_C,\\
	r_i &\doteq I_{\mathcal{C}(i)} r,\\
	p_i &\doteq I_{\mathcal{C}(i)} p.
\end{align*} Let $\Lambda \doteq \sum_{i \in \mathcal{M}} I_{\mathcal{C}(i)}^{\top} I_{\mathcal{C}(i)}$ and $\Lambda_i \doteq  I_{\mathcal{C}(i)} \Lambda I_{\mathcal{C}(i)}^{\top}$. We rewrite (\ref{eq:CGalpha}) as $\alpha^n = \frac{\eta^n}{\sigma^n}$ where $\eta^n \doteq r^{n\top}r^n$ and $\sigma^n \doteq p^{n\top}Sp^n$. This can be decomposed into 
\begin{align*}
\eta^n &= \sum_{i \in \mathcal{M}} \eta^n_i, \quad \eta^n_i \doteq r_i^{n\top} \Lambda_i^{-1} r_i^n,\\
\sigma^n &= \sum_{i \in \mathcal{M}} \sigma_i^n, \quad \sigma_i^n \doteq p^{n\top}_i \hat{S}_i p^n_i,
\end{align*} where $\hat{S}_i \doteq I_{\mathcal{C}(i)} S_i I_{\mathcal{C}(i)}^\top $. The computation of $\alpha^n$ therefore requires the local computation of $\eta_i$ and $\sigma_i$ for each agent and two subsequent scalar global sums. 
We next decentralize the updates of the Lagrange multiplier $\lambda_C$ and the step direction $p$. We multiply (\ref{eq:CGlambda}) and (\ref{eq:CGp}) by $I_{\mathcal{C}(i)}$ from the left and get
\begin{align*}
\lambda_{C,i}^{n+1} &= \lambda_{C,i}^{n} + \frac{\eta^n}{\sigma^n} p^n_i,\\
p^{n+1}_i &= r_i^{n+1} + \beta^n p_i^n,
\end{align*} where $\beta^n = \frac{\eta^{n+1}}{\eta^{n}}$ as before. 
Note that both updates can be performed locally. In the last step, we decompose the residual update
\begin{align*}
r^{n+1}_i = r^n_i - \frac{\eta^n}{\sigma^n}\sum_{j\in \mathcal{M}(i) \cup i} I_{ij} \hat{S}_j   p^n_j,
\end{align*} where $I_{ij}\doteq I_{\mathcal{C}(i)} I_{C(j)}^\top.$   This requires local communication among neighboring agents. We have now decomposed all CG updates and can summarize DCG in Algorithm \ref{alg:d-CG}. The initialization of the residual and step direction must satisfy $r^0=p^0=s -S  \lambda_C^0$ where $\lambda_C^0$ can be chosen. This can be decomposed into
\begin{equation*}
r_i^0 = p_i^0 = \sum_{j \in \mathcal{M}(i) \cup i} I_{ij} I_{C(j)}s_j - \sum_{j \in \mathcal{M}(i) \cup i} I_{ij}\hat{S}_j \lambda_{C,j}^0 
\end{equation*} and hence requires neighbor-to-neighbor communication.

\begin{algorithm}[t]
	\caption{DCG}
	\textbf{Initialization: $ \lambda_C^0$ and $r^0=p^0=s -S  \lambda_C^0$}\\
	\textbf{Repeat until} $\lVert r_i^n \rVert < \varepsilon \quad \forall i \in \mathcal{M}$ \textbf{:}
	\begin{enumerate}\setlength\itemsep{-0.1cm}
		\item$\eta_i^n = r_i^{n\top} \Lambda_i^{-1} r_i^n, \quad \sigma^n_i =  p^{n\top}_i \hat{S}_i p^n_i$\\
		\item$\eta^n = \sum_{i \in \mathcal{M}} \eta^n_i, \quad \sigma^{n} = \sum_{i \in \mathcal{M}} \sigma_i^{n}$\\ 
		\item$\lambda_{C,i}^{n+1} = \lambda_{C,i}^{n} + \frac{\eta^n}{\sigma^n} p^n_i$\\	
		\item$r^{n+1}_i = r_i^n - \frac{\eta^{n}}{\sigma^{n}} \sum_{j \in \mathcal{M}_i \cup i} I_{ij} \hat{S}_j   p^n_j$\\
		\item$p^{n+1}_i = r_i^{n+1} + \frac{\eta^{n+1}}{\eta^{n}} p_i^n$\\
	\end{enumerate} \label{alg:d-CG}
\end{algorithm}

\subsection{Communication analysis}
We rely on three measures to analyse the communication requirements of ASM/DCG: the number of floats that are sent to and from a central coordinator (global floats), the number of booleans that are sent to and from a central coordinator (global booleans), and the number of floats that are sent on a neighbor-to-neighbor basis (local floats). In each DCG iteration, each agent sends $\eta_i^n$ and $\sigma_i^n$ to a central coordinator. The coordinator then computes the sums $\sigma^n$ and $\eta^n$ and returns them to each agent. This totals in $4M$ global floats per DCG iteration. In addition, each agent sends a convergence flag to the coordinator and the coordinator returns a global convergence flag to each agent, which gives $2M$ global booleans per DCG iteration. In step 4. of the DCG iteration, agent $j$ sends those elements of $\hat{S}_j^n p_j^n$ to neighbor $i$, that correspond to non-zero elements in $I_{ij}$. Hence, agents $i$ and $j$ exchange one float per coupled variable. This results in $2n_c$ local floats per DCG iteration. Within each ASM iteration, either the smallest Lagrange multiplier or the largest feasible step size have to be determined and convergence is checked. This requires $2M$ global floats and $2M$ global booleans. The communication footprint for the ASM, DCG and the ADMM variant explained in Appendix \ref{app:ADMM} is summarized in Table \ref{tab:commAnalysis}.

\begin{table}[b]
	\begin{center}
		\caption{Communication per iteration.}\label{tab:commAnalysis}
		\begin{tabular}{p{0.9cm} c c c }	
			& global fl. & global bo. & local fl.\\
			\hline
			DCG & $4M$ & $2M$ & $2n_c$\\
			ASM & $2M$ & $2M$ & $0$\\
			ADMM & $0$ & $2M$ & $2n_c$
		\end{tabular}
	\end{center}
\end{table}

\section{Numerical Results}
\label{sec:Results}
We compare the presented ASM with an ADMM  implementation (cf. Appendix \ref{app:ADMM}) on a chain of masses system \citep{Conte2016}. The baseline example consists of 10 masses with $m=\SI{1}{\kg}$ that are coupled by springs with stiffness $k=\SI{3}{\N/ \m}$ and dampers with coefficient $d = \SI{3}{\N\s/\m}.$ Each mass is actuated by a force $\SI{-1}{\N} \leq u_i \leq \SI{1}{\N}.$ We do not consider constraints on the states $x_i(t)=\begin{bmatrix} y_{i}(t) & v_{i}(t) \end{bmatrix}^\top$, where $y_{i}(t)$ and $v_{i}(t)$ denote the position and velocity of mass $i$. The equations of motion are discretized using the Euler forward discretization with step size $T=\SI{0.2}{\s}.$ We choose the parameters of the MPC controller as $N=12$, $Q_i = \text{diag}([10,10])$, $P_i = 0$ and $R=1$. The resulting OCP (\ref{ocp:distributed}) has the dimension $n_z=812$ with 260 equality constraints for the dynamics, 240 input constraints and 432 coupling constraints. We choose the following tolerances for the stopping criteria: $\lVert r_i \rVert_\infty < 10^{-7}$ for DCG and $\lVert \Delta z_i \rVert_{\infty} < 10^{-6}$ for ASM. We choose $\lambda_{C,i}^0=0$ to initialize DCG and warm-start the ASM with the optimal active set from the previous MPC iteration.
For ADMM, the stopping criteria $\lVert C_i^c(z_i-\bar{z}_i) \rVert_\infty \leq \varepsilon_r \text{min}\{\text{max}\{\lVert C_i^c z_i \rVert_\infty, \lVert C_i^c \bar{z}_i \rVert_\infty\},1\}$ and $\lVert \rho C_i^c(z_i^+ - z_i) \rVert_\infty \leq \varepsilon_d \text{min}\{\lVert \lambda_{C,i}\rVert_\infty,1\}$ are used. Two different tolerance levels are chosen: 

\begin{itemize}\setlength\itemsep{-0.4cm}
	\item ADMM1: $\varepsilon_r=10^{-6}$ and $\varepsilon_d=10^{-3}$\\
	\item ADMM2: $\varepsilon_r = 10^{-4}$ and $\varepsilon_d = 10^{-2}.$
\end{itemize}

\begin{figure}[t]
	\centering
	\pgfdeclarelayer{background}
	\pgfdeclarelayer{foreground}
	\pgfsetlayers{background,main,foreground}
	\input{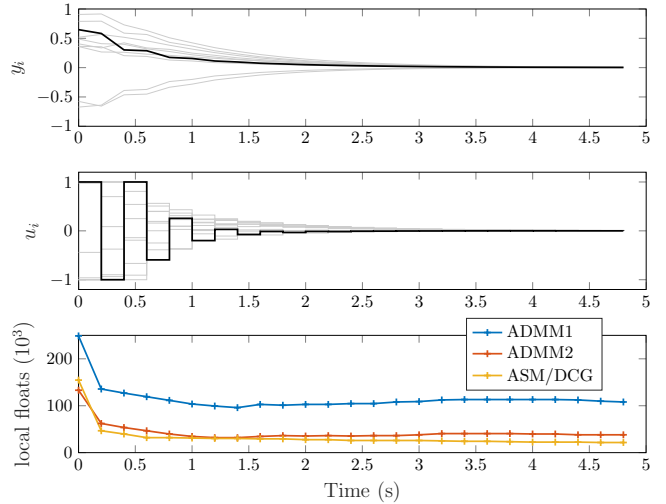}
	\label{fig:asmVsAdmm}
	\caption{Closed-loop trajectories and locally communicated floats per sampling instance.}
\end{figure}

Five case studies that each simulate the closed-loop behavior for 25 MPC iterations starting at 30 random initial positions are conducted. For the baseline case study, the initial position and velocity for each agent are chosen from a uniform distribution in the invertvals $\SI{-1}{\meter} \leq y_{i,0} \leq \SI{1}{\meter}$ and $\SI{-0.5}{\meter / \second} \leq v_{i,0} \leq \SI{0.5}{\meter / \second}.$ Figure \ref{fig:asmVsAdmm} shows simulation results with the closed-loop trajectories of the agents and the number of floats per sample interval that is communicated locally from neighbor-to-neighbor in the network. The output and input trajectories of the fifth agent are displayed in black for better visualization and the trajectories of the remaining agents are shown in grey.
Table \ref{tab:baselineFloats} shows the communication in terms of floats and booleans sent to and from a central coordinator (global floats) as well as floats sent on a neighbor-to-neighbor basis (local floats). We note that this table and all further tables do not include data on the first MPC iteration for each initial position as these iterations cannot be warm-started. We further note that ASM/DCG achieves an accuracy of $10^{-7}$ for the closed loop state trajectories compared to the centralized MPC scheme whereas ADMM1 and ADMM2 achieve only accuracy levels of $10^{-5}$ and $10^{-4}$, respectively. Table \ref{tab:baselineIter} shows the number of iterations required by each algorithm per sample interval. For ASM/DCG, we count outer iterations spent in Algorithm 1 (ASM) and inner iterations spent in Algorithm \ref{alg:d-CG} (DCG) separately. The iterations split into a first part to obtain $z_i^0$ and into a second part to update the active set. For the baseline case, most ASM/DCG iterations are spent finding a feasible $z_i^0$. Future work may improve the initialization procedure (cf. Remark \ref{rema:initASM}).
Four additional case studies were conducted in addition to the baseline case. A single parameter compared to the baseline case has been changed for each of the additional cases and the results are given in Table \ref{tab:furtherStudies}. The iterations reported for ASM/DCG are inner (i.e. DCG) iterations. The results show that ASM/DCG requires less iterations than ADMM1 for all analyzed scenarios.

\begin{table}
\begin{center}
	\caption{Communication footprint (baseline case).}\label{tab:baselineFloats}
\begin{tabular}{p{1.5cm} p{0.7cm} p{0.7cm} p{0.7cm} p{0.7cm} p{0.7cm} c }	
	& \multicolumn{2}{c}{global floats} & \multicolumn{2}{c}{global bool.} & \multicolumn{2}{c}{local floats}\\
	\hline
	& \centering mean & \centering max & \centering mean & \centering max & \centering mean & max\\
	ASM/DCG & \centering 1.3k & \centering 3.9k & \centering 0.7k & \centering 2.1k & \centering 27k & 88k\\
	ADMM1 & \centering 0 & \centering 0 & \centering 2.4k & \centering 3.7k & 102k & 160k\\
	ADMM2 & \centering 0 & \centering 0 & \centering 0.8k & \centering 1.6k & \centering 35k & 68k
\end{tabular}
\end{center}
\end{table}

\begin{table}
	\begin{center}
		\caption{Iterations  (baseline case).}\label{tab:baselineIter}
		\begin{tabular}{p{1.5cm}  p{0.7cm} p{0.7cm} p{0.7cm} p{0.7cm} p{0.7cm} c }	
			& \multicolumn{2}{c}{total iter.}& \multicolumn{2}{c}{feas. guess} & \multicolumn{2}{c}{AS updating}\\
			\hline
			& \centering mean & \centering max & \centering mean & \centering max & \centering mean & max\\
			DCG & \centering 30& \centering 98 & \centering 27 & \centering 97 & \centering 1 & 1\\
			ASM & \centering 1 & \centering 1 & \centering - & \centering - & \centering 1 & 1\\
			ADMM1 & \centering 117 & \centering 185 & \centering - & \centering - & \centering - &-\\
			ADMM2 & \centering 41 & \centering 78 & \centering - & \centering - & \centering - &-\\
		\end{tabular}
	\end{center}
\end{table}

\begin{table}
	\begin{center}
		\caption{Iterations and communication footprint (modified problem parameters).}\label{tab:furtherStudies}
		\begin{tabular}{p{1.5cm}  p{0.7cm} p{0.7cm} p{0.7cm} p{0.7cm} p{0.7cm} c }	
			& \multicolumn{2}{c}{total iter.}& \multicolumn{2}{c}{feas. guess} & \multicolumn{2}{c}{local floats}\\
			\hline
			& \centering mean & \centering max & \centering mean & \centering max & \centering mean & max\\
			\tiny{$\lvert v_{i,0} \rvert \leq 2$}& & & & &\\
			ASM/DCG & \centering 37 & \centering 283 & \centering 29 & \centering 108 & \centering 34k &251k \\
			ADMM1 & \centering 133 & \centering 303 & \centering - & \centering - & \centering 109k & 262k\\
			\hline
			\tiny{$N=5$} & & & & &\\
			ASM/DCG & \centering 30 & \centering 137 & \centering 29 & \centering 108 & \centering 12k &52k \\
			ADMM1 & \centering 154 & \centering 226 & \centering - & \centering - & \centering 56k & 82k\\
			\hline
			\tiny{$5$ Masses} & & & & &\\
			ASM/DCG & \centering 26 & \centering 68 & \centering 25 & \centering 67 & \centering 11k &28k \\
			ADMM1 & \centering 106 & \centering 171 & \centering - & \centering - & \centering 41k & 66k\\
			\hline
			\tiny{$20$ Masses} & & & & &\\
			ASM/DCG & \centering 32 & \centering 160 & \centering 29 & \centering 153 & \centering 61k & 301k \\
			ADMM1 & \centering 127 & \centering 237 & \centering - & \centering - & \centering 231k & 433k\\
		\end{tabular}
	\end{center}
\end{table}

\section{Summary and Outlook}
This paper has proposed a novel distributed active set method that can be used for the distributed MPC of linear systems. The method combines a usual active set approach with a recently proposed decentralized conjugate gradient method to solve the arising equality constrained QPs. It thus is a distributed method. Our numerical case study shows a competitive communication footprint compared to an ADMM variant in terms of number of floats communicated on a neighbor-to-neighbor basis. Future work will include the extension to state and terminal constraints. In particular, the proposed scheme can be terminated early for MPC schemes involving terminal constraints and still guarantee stability because the scheme ensures primal feasibility. This aspect and improved initialization procedures may also reduce the communication requirement further.

{\footnotesize
\bibliography{ifacconf}}             % bib file to produce the bibliography
                                                     % with bibtex (preferred)
                                                   
\appendix
\section{ADMM} \label{app:ADMM}
We briefly comment on the ADMM variant that we use here for the comparison with ASM/DCG. We refer the reader to \citep{Boyd2011} for an extensive overview of ADMM and to \citep{Rostami2017} for the application of ADMM to distributed MPC. To apply ADMM, we introduce the trajectory average $\bar{\boldsymbol{x}}_i \in \mathbb{R}^{Nn_i}$ and the decision variable $\bar{z}_i^\top \doteq [\bar{\boldsymbol{x}}_i^\top, x_i^{N\top}, \boldsymbol{u}_i^\top, [\bar{\boldsymbol{x}}_j^\top]_{j \in \mathcal{M}_i^{\text{in}}}] \in \mathbb{R}^{n_{zi}}$, which satisfy the coupling constraints (\ref{eq:cCon}) by design. ADMM alternately updates $z_i$, which satisfies (\ref{eq:eqCon}) and (\ref{eq:IneqCon}), and $\bar{z}_i$ until $z_i$ satisfies (\ref{eq:cCon}) to a chosen accuracy. The  $\bar{z}_i$ update requires local communication and convergence flags are sent to a global coordinator as for ASM/DCG.

\begin{algorithm}[h!]
	\caption{ADMM}
	\textbf{Initialization: $ \lambda_{C,i}^0$ and $\bar{z}_i^0$}\\
	\textbf{Repeat until convergence:}
	\begin{enumerate}\setlength\itemsep{-0.05cm}
		\item $z_i^{n+1} = \displaystyle\argmin_{z_i} z_i^{\top}H_i z_i + \lambda_{C,i}^{\top} C_i^c z_i + \frac{\rho}{2}\lVert C_i^c(z_i-\bar{z}_i^n)\rVert_2^2$\\
		\hspace*{1.4cm} s.t. $\;\;C_i^{\mathcal{E}}z_i = b_i^{\mathcal{E}}, \;\; C_i^{\mathcal{I}}z_i \leq b_i^{\mathcal{I}}$
		\item Receive $[\boldsymbol{v_{ij}}]_{j \in \mathcal{M}_i^{\text{out}}}.$
		
		\item Compute $\bar{\boldsymbol{x}}_i = \sum_{j \in \mathcal{M}_i^{\text{out}}} (\boldsymbol{x_i}+\boldsymbol{v_{ij}})/(2 \lvert \mathcal{M}_i^{\text{out}} \rvert) $
		
		\item Receive $[\bar{\boldsymbol{x}}_j]_{j \in \mathcal{M}_i^{\text{in}}} $ and form $\bar{z}_i$.
		
		\item$\lambda_{C,i}^{n+1} =\lambda_{C,i}^n + \rho C_i^c (z_i^{n+1}-\bar{z}_i^{n+1})$\\
	\end{enumerate} \label{alg:ADMM}
\end{algorithm}

\end{document}